\documentclass[11pt,a4paper]{article}

\usepackage{a4}

\usepackage{amssymb}

\usepackage{amsmath}

\usepackage{euscript}

\author{Leonard Todjihounde}

\title{Mean-value property on manifolds with minimal horospheres}

\begin{document}

\date{}

\maketitle

\begin{abstract}

\noindent Let $(M,g)$ be a non-compact and complete Riemannian
manifold with mini\-mal horospheres and infinite injectivity radius.
 We prove that bounded functions on $(M,g)$ satisfying the mean-value
property are constant. \\We extend thus a result of the authors in
[6] where they proved a similar result for bounded harmonic
functions on harmonic manifolds with minimal horospheres. 

\end{abstract}

\section{Introduction}

\noindent Let $(M,g)$ be a non-compact and complete Riemannian
manifold. \\ A function $u$ defined on $(M,g)$ is said to have the
mean-value property if:
\[\forall \;r > 0 \;\mbox{and}\;\forall \;p \in M \;,\; u(p) =
\frac{1}{V(p,r)}\int_{B(p,r)}u(q)\;d\mu (q)\;,\] where $d\mu$
denotes the Riemannian volume element and $V(p,r)$ the volume of the
closed ball $B (p,r)$ of centre $p$ and radius $r$. \\ Well-known
examples of functions satisfying the mean-value
property are harmonic functions on harmonic manifolds (see [9]). \\
In [6] the authors proved that on non-compact harmonic
 manifolds with minimal horospheres, bounded harmonic functions are
 constant.\\ One of the major arguments to obtain this result is
 the fact that on harmonic manifolds, harmonic functions possess
 the mean-value property. It thus seems natural to raise the same
 question by considering the class of functions satis\-fying the
 mean-value property and defined on manifolds not necessarily harmonic. \\
 Some analogous of Liouville type results for functions
 satisfying the mean-value property have been proved by several
 authors. For example the authors in [8] proved that on certain
 kinds of homogeneous spaces, the only $L^p$-function possessing
 the mean-value property is the zero function. For similar results
 and related works see also [1, 2, 3, 4, 5, 10] and the references
 therein. \\
 Our aim is to extend the Liouville type result proved in [6] on bounded
 functions satysfying the mean-value property and defined on non-compact manifolds with
 minimal horospheres and infinite injectivity radius. \\
 We refer to [6] and [7] for information and
 details on the minimality's condition of horospheres in a
 non-compact manifold. \\
\noindent For a real number $r > 0$, we consider the stability
vector field $H(.,r)$ defined by:
\[H (p,r) =: \int_{B(p,r)}\exp_p^{-1}(q)\;d\mu (q) \;,\;\forall \; p \in M \;,\]
where $\exp_p^{-1}$ denotes the inverse of the exponential map.\\
Let us note that the vanishing of the stability vector field for any
radius $r > 0$ means that any geodesic ball in $(M,g)$ has its
Riemannian center of mass (or center of gravity) at the centre of
the ball. This is the case for examples for harmonic manifolds,
d'Atri spaces or compact locally symmetric spaces. \\
In the next section we give a result relating the gradient of the
volume function and the derivative of the stability vector field
that we use in the third section to prove that on non-compact
manifolds with minimal horospheres and infinite injectivity radius,
bounded functions having the mean-value property are constant.

\section{Volume functions and stability vector fields}

\noindent Let $V : (p,r) \in M \times \;]0\;,\;+\infty[\;
\longmapsto V(p,r)$ be the function associating to each pair $(p,r)
\in M \times \; ]0\;,\;+\infty[$ the volume $V(p,r)$ of the ball
$B(p,r)$. \\The volume function $V$ and the stability vector
field $H$ are related by the following differential equation: \\ \\
\noindent {\bf Lemma 2.1} \\
Let $\nabla$ denotes the gradient operator on $(M,g)$.\\
For any $r > 0$ and $p \in M$, it holds:
\[\nabla V(p,r) - \frac{1}{r}\frac{\partial}{\partial r}H (p,r)  = 0\;.\]
\\
\noindent {\bf Proof}  \\
For $X \in T_pM$, it holds
\[\nabla_X V(p,r) = \int_{S(p,r)} <\eta (q)\;,\;X(q)>\;d\sigma (q)
\;,\] where $d\sigma$
 denotes the Riemannian measure induced on the sphere $S(p,r)$ with centre
$p$ and radius $r$, $\eta (q)$ is the outward unit normal at $q$,
and $X(q)$ is the parallel transport of $X$ from $p$ to $q$. \\
By the Gauss Lemma,
\[<\eta (q)\;,\;X(q)>\; = \;<(d\exp_p^{-1}) \eta (q)\;,\;X> \;.\]
Otherwise
\[(d\exp_p^{-1}) \eta (q) = r^{-1}\exp_p^{-1}q \;.\]
It follows then
\begin{eqnarray*}
\nabla_X V(p,r) &=& \int_{S(p,r)}r^{-1}\exp_p^{-1}q \;d\sigma (q)
\\ &=& r^{-1}\int_{S(p,r)}\exp_p^{-1}q \;d\sigma (q) \\
&=& r^{-1}\frac{\partial}{\partial r}\left(
\int_{B(p,r)}\exp_p^{-1}q \;d\mu (q) \right) \\
&=& r^{-1}\frac{\partial}{\partial r}H(p,r) \;.
\end{eqnarray*}
Hence the result.  \hfill$\square$

\section{A derivative formula}

Let $p \in M$ and $X$ a unit vector in $T_pM$. We consider as in [6]
the function:
\begin{eqnarray*}
\theta_X : M - \{p\} &\longrightarrow& {\mathbb R} \\
q &\longmapsto& \theta_X(q) =: \angle_p (X ,
\overset{.}{\gamma}_q(0)) \;,
\end{eqnarray*}
where $\gamma_q$ denote the geodesic defined by $\gamma_q (t) =
\exp_p (t\exp_p^{-1}q) \;,\;\forall \;t \in [0\;,\;1]$, and
$\angle_p (X , \overset{.}{\gamma}_q(0))$ the angle at $p$ between
the vectors $X$ and $\overset{.}{\gamma}_q(0)$. \\For the geodesic
$c$ with $c (0) = p$ and $\overset{.}{c}(0) = X$, let $P_t$ be the
parallel transport along $c$ and $f_t$ the one parameter family of
diffeomorphisms of $M$ given by
$f_t = \exp_{c(t)}\circ P_t \circ \exp_p^{-1}$. \\
Let $u$ be a differentiable function on $M$ possessing the
mean-value property. \\
It holds: \\ \\
\noindent {\bf Proposition 3.1}\\
For any real number $r > 0$,
\[ X u (p) = \frac{1}{V(p,r)}\int_{S(p,r)} u \cos \theta_X d\sigma
 - \frac{1}{r}\frac{u(p)}{V(p,r)}<\frac{\partial}{\partial r}H(p,r)\;,\;X> \;.\]
\\
\noindent {\bf Proof} \\
Since the function $u$ possesses the mean-value property, we have:
\[u(c(t)) = \frac{1}{V(c(t),r)}\int_{B(c(t),r)}u(q)\;d\mu (q)\;.\]
And then
\begin{eqnarray*}
X.u (p) &=& \frac{d}{dt}u(c(t))_{|t=0} \\
&=& \frac{d}{dt}\left(\frac{1}{V(c(t),r)}\int_{B(c(t),r)}u\;d\mu
\right)_{|t=0} \\
&=& - \frac{1}{V(p,r)^2} <\nabla V(p,r)\;,\;X> \int_{B(p,r)} u
\;d\mu \\ && {} + \frac{1}{V(p,r)} \frac{d}{dt}
\left(\int_{B(c(t),r)}u\;d\mu \right)_{|t=0} \;\;\;\;\;\;(i) \;.
\end{eqnarray*}
From Lemma 2.1,
\[\nabla V(p,r)  = \frac{1}{r}\frac{\partial}{\partial r}H
(p,r)\;.\] Thus
\begin{eqnarray*}
\frac{1}{V(p,r)^2} <\nabla V(p,r),X> \int_{B(p,r)} u \;d\mu &=&
\frac{1}{V(p,r)^2} <\frac{1}{r}\frac{\partial}{\partial r}H
(p,r),X> \int_{B(p,r)} u \;d\mu \\
&=& \frac{1}{r}\frac{u(p)}{V(p,r)}<\frac{\partial}{\partial
r}H(p,r),X> \;\;\;\;(ii), \\
&& {} \mbox{since}\;\;u(p) = \frac{1}{V(p,r)}\int_{B(p,r)}u(q)\;d\mu
(q)\;.
\end{eqnarray*}
By Theorem 2.1 in [6] we have:
\begin{eqnarray*}
\frac{d}{dt}\left(\int_{B(c(t),r)}u\;d\mu\right)_{|t=0} &=&
\frac{d}{dt}\left(\int_{B(p,r)} f_t^*(u\;d\mu) \right)_{|t=0} \\
&=& \int_{B(p,r)}\frac{d}{dt}(f_t^*(u\;d\mu))_{|t=0} \\
&=& \int_{S(p,r)}u \cos \theta_X \;d\sigma  \;\;\;\; (iii)\;.
\end{eqnarray*}
By replacing $(ii)$ and $(iii)$ in the relation $(i)$ we obtain the
result. \hfill$\square$
\\ \\
\noindent By using the derivative formula given in Proposition 3.1, we get: \\ \\
\noindent {\bf Theorem 3.1} \\
Let $(M,g)$ be a non-compact and complete Riemannian manifold
with minimal horospheres and infinite injectivity radius. \\
Any bounded function on $(M,g)$ satisfying the mean-value property
is constant.
\\ \\
\noindent {\bf Proof} \\
Let $u$ be a bounded function on $(M,g)$ satisfying the mean-value
property.
\\ By Proposition 3.1,
\[|X u (p)| \le \alpha \frac{A(p,r)}{V(p,r)}
 + \frac{\alpha}{V(p,r)} \|\frac{1}{r}\frac{\partial}{\partial
 r}H(p,r)\|\;,\;\forall \; p \in M \;\mbox{and}\;r > 0 \;,\]
 where $A(p,r)$ is the area of the sphere $S(p,r)$, and
 $\alpha \ge 0$ is such that $|u| \le \alpha $. \\
 But:
\begin{eqnarray*}
\|\frac{1}{r}\frac{\partial}{\partial
 r}H(p,r)\| &=& \|\frac{1}{r} \frac{\partial}{\partial
 r} \int_{B(p,r)}\exp_p^{-1}q\;d\mu (q)\| \\
 &=& \|\frac{1}{r} \int_{S(p,r)}\exp_p^{-1}q\;d\sigma (q)\| \\
 &\le& \frac{1}{r} \int_{S(p,r)}\|\exp_p^{-1}q \|\;d\sigma (q)\\
 &=& A(p,r)\;,\;\mbox{since}\;
 \|\exp_p^{-1}q \| = r \;,\; \forall \;q \in S(p,r)\;.
\end{eqnarray*}
 Thus we get:
\[|X u (p)| \le 2\alpha \frac{A(p,r)}{V(p,r)}\;.\]
 Due to the minimality of horospheres (see [6] for details)
\[\underset{r \to + \infty}{\lim}\frac{A(p,r)}{V(p,r)} = K_\infty = 0 \;.\]
By taking the limit of the previous inequality as $r \to \infty$, it
follows then:
\[|X u (p)| = 0\;,\;\;\mbox{for any}\;\; p \in M \;\mbox{and any unit vector}
\; X \in T_pM \;.\] Hence $u$ is a constant function.
\hfill$\square$  \\ \\
\noindent {\bf Remark}: From the proof of Theorem 3.1 it is easy to
see that the same result can be obtained by assuming the horospheres
with bounded mean-value and not necessarily minimal.

\end{document}